\documentstyle[art10]{article}
\begin{document}

\newtheorem{lemma}{Lemma}
\newtheorem{theorem}{Theorem}
\newtheorem{definition}{Definition}
\newtheorem{corollary}{Corollary}
\newtheorem{proposition}{Proposition}

\def\R{{\bf R}}
\def\C{{\bf C}}
\def\Z{{\bf Z}}
\def\Q{{\bf Q}}

\begin{center}
{\bf On existence of nonformal simply connected symplectic
manifolds}
\footnote{This work is supported by RFBR, grants
96-01-00182a and 98-01-00749).}
\end{center}

\begin{center}
{\bf Ivan K. Babenko and Iskander A. Taimanov}
\end{center}

A smooth manifold is called symplectic if there is
a nondegenerate closed 2-form $\omega$ on it.  By the
Darboux theorem near every point this form is
reduced to the form
$\omega = \sum_{j=1}^n d x^j \wedge d x^{j+n}$ with
$\{x^j\}$ local coordinates on the manifold. This implies that
the dimension of a symplectic manifold
$M$ is even and there is an almost complex structure on it:
take a Riemannian metric $(\cdot ,\cdot )$ on $M$ and define
by $\omega(u,Jv) = (u,v)$ an automorphism $J:TM \rightarrow TM$
defining an almost complex structure.
In late 60s Gromov showed that if the cohomology class of
$\omega$ on a compact symplectic manifold $(M,\omega)$ is integer,
$[\omega] \in H^2(M;\Z)$, then for sufficiently large
$N$ there exists an embedding $f:M
\rightarrow \C P^N$ such that $f^*(\omega_0) = \omega$ with
$\omega_0$ the Hodge form on $\C P^N$
(\cite{Gromov}). As it was shown by Tischler such embedding
already exists for $N = 2n+1$ with $\dim M = 2n$ (\cite{Tischler}).

Since K\"ahler manifolds are symplectic (take for $\omega$
the K\"ahler form), a conjecture that all
compact symplectic manifolds are exhausted by K\"ahler manifolds
up to diffeomorphisms appears.

For non-simply-connected manifolds this had been disproved by
Thurston who had proposed for a simplest example
$\tilde{N} = \R^3/\Gamma  \times S^1$ with $\Gamma$
a uniform lattice of uppertriangular integer
$3\times 3$-matrices in the Heisenberg group
(\cite{Thurston}).  It appeared that before this example
had been found by Kodaira. Notice that the symplectic form on
$\tilde{N}$ is integer.

An example of simply connected symplectic manifold
nonhomeomorphic to a K\"ahler manifold was constructed by
McDuff who had embedded the Kodaira--Thurston manifold
$\tilde{N}$ into $\C P^5$ and had constructed a manifold
$X$ as a symplectic blowing up of $\C P^5$ at $\tilde{N}$
(\cite{McDuff}). The latter procedure consists in
fiberwise blowing up of the normal bundle to
$\tilde{N}$ in $\C P^5$ (since
$\tilde{N}$ is a symplectic submanifold, fibers of the
normal bundle are endowed with the canonical almost
complex structure).

Later Gompf had even constructed four-dimensional simply
connected symplectic manifolds which are not homeomorphic to
K\"ahler manifolds (\cite{Gompf}).

An important property of K\"ahler manifolds established in
\cite{DGMS} is their formality. A polyhedron
$M$ is called formal if there exists a homomorphism of
skew-commutative graded differential algebras
$({\cal M}(M),d) \rightarrow (H^{\ast}(M;\Q),d)$ such that
it induces an isomorphism of cohomologies. Here
$(H^{\ast}(M),d)$ is a ring of rational cohomologies
endowed with a zero differential $d \equiv 0$ and $({\cal
M}(M),d)$ is the minimal model of $M$ (for the definition of
minimal models of simply connected and nilpotent spaces see, for
instance, \cite{DGMS}).

Being poor with in methods of constructing a collection of
examples of symplectic manifolds
generates many papers on finding symplectic non-K\"ahler
non-simply-connected manifolds (see, for instance,
\cite{CFG,BG,LO}) and therewith the formality criterion was
used in \cite{BG,LO}. A problem on existence of nonformal
simply connected symplectic manifolds remained open.
The most complete collection of results on formality of
symplectic manifolds is given in \cite{LO}.

In general, as we show, a simply connected symplectic manifold
is nonformal. Therefore despite a strong
differential-geometric structure which is a symplectic form
the topology of a simply connected symplectic manifold can
be substantially more complicated than topology of K\"ahler
manifolds.

We define $X_k$ as a symplectic blowing up of $\C P^k$
at embedded the Kodaira--Thurston manifold $\tilde{N}$.
In these notations the McDuff example is $X_5$.  All manifolds
$X_k$ are simply connected for $k \geq 5$.

The following theorem holds.

{\bf Theorem.}
{\sl For $k \geq 6$ the manifold $X_k$ is nonformal.}

The proof is as follows.

Let $\pi:X_k \rightarrow \C P^k$ be the projection
whose restriction onto the complement to the
preimage of $\tilde{N}$ is a diffeomorphism onto the image and
let $\tilde{V}$ be a tubular neighborhood of $\pi^{-1}(\tilde{N})$.
The exact cohomology sequence of the pair
$(X_k,\tilde{V})$ contains the following fragment
\begin{equation}
0 \leftarrow
H^4(X_k,\tilde{V}) \stackrel{\partial}{\leftarrow}
H^3(\tilde{V}) \approx H^3(\tilde{N}) \oplus \mbox{Ker}\,
\partial \leftarrow H^3(X_k) \leftarrow 0.
\label{1}
\end{equation}

A space $\tilde{V}$ is diffeomorphic to a $\C
P^{k-3}$-bundle over $\tilde{N}$ and $H^{\ast}(\tilde{V})$
is isomorphic to a skewed tensor product
$H^{\ast}(\tilde{N}) \hat{\otimes} H^{\ast}(\C P^{k-3})$.
Denote by $v$ some generator of
$H^1(\tilde{N})$, denote by $b$ the generator of
$H^2(\C P^{k-3})$, denote by $a \in H^2(X_k)$ the cohomology
class of the symplectic form, and denote by
$u \in H^3(X_k)$ an element which is mapped into $b \wedge v$
(\ref{1}). It is shown that $u \wedge a =
0$ and, since $u$ and $a$ are spheric cycles, there
exists a generator $\hat{z}$, of the minimal model
${\cal M}(X_k)$ of the space $X_k$, such that
$d\hat{z} = \hat{u} \wedge \hat{a}$ where $\hat{u}$ and
$\hat{a}$ are mapped into
$u$ and $a$ by the Hurewicz homomorphism.

It is shown by straightforward computations that
$\hat{z} \wedge \hat{u} \in {\cal M}(X_k)$
generates a nontrivial seven-dimensional cohomology cycle
and therefore $X_k$ is nonformal.
For these computations it is essentially that
$k \geq 6$.

The detailed proof will be published elsewhere.

\end{document}